\documentclass[preprint,12pt,notitlepage]{article}

\usepackage{amsmath,amssymb}
\usepackage{amsthm}

\usepackage{enumerate} %
\usepackage{tikz-cd}
\usepackage{tikz}
\usepackage{pgf}
\usepackage{xifthen}

\tikzset{mk@vertex/.style={font=\scriptsize}}

\newtest{\IfVal}[4]{{#1=#2} \and {#3=#4}}
\newcommand{\dywiz}{\kern0sp\discretionary{-}{-}{-}\penalty10000\hskip0sp\relax}

\newcommand{\add}{\operatorname{add}} 
\newcommand{\ann}{\operatorname{ann}} 
\newcommand{\D}{\operatorname{D}} 
\newcommand{\modd}{\operatorname{mod}} %
\newcommand{\ind}{\operatorname{ind}} %
\newcommand{\dtr}{\operatorname{DTr}}
\newcommand{\trd}{\operatorname{TrD}}
\newcommand{\End}{\operatorname{End}} 
\newcommand{\Ext}{\operatorname{Ext}} 
\newcommand{\Hom}{\operatorname{Hom}} 
\newcommand{\rad}{\operatorname{rad}} 
\newcommand{\rk}{\operatorname{rk}} 
\newcommand{\soc}{\operatorname{soc}} 
\renewcommand*\top{\operatorname{top}} 

\def\adm#1{$(\mathop{ad} #1)$} %
\providecommand{\abs}[1]{\lvert#1\rvert}
\newcommand{\bbA}{\mathbb{A}} %
\newcommand{\bbX}{\mathbb{X}} %
\newcommand{\bbQ}{\mathbb{Q}} %
\newcommand{\bbZ}{\mathbb{Z}} %
\newcommand{\hB}{\widehat{B}} %
\newcommand{\mc}{\mathcal} %
\newcommand{\wt}[1]{\widetilde{\mathbb{#1}}} %

\newcommand{\LinGraph}[2]
{
\begin{tikzpicture}
\foreach \x in {1,...,#1} { \node (e\x) at (1.2*\x,0) {$\bullet$}; }%

\foreach \x in {2,...,#1} {%
\pgfmathparse{int(\x-1)}\let\y\pgfmathresult;%
\draw[-] (e\x)--(e\y) node[mk@vertex, pos=.5,sloped,above] {\foreach
\i/\j/\valA/\valB in #2
{\ifthenelse{\IfVal{\j}{\x}{\i}{\y}}{$(\valA,\valB)$}{}}}; }
\end{tikzpicture}
}

\newcommand{\DGraph}[2]
{
\begin{tikzpicture}[scale=0.8]
\pgfmathparse{int(#1)}\let\c\pgfmathresult;%
\pgfmathparse{int(#1-1)}\let\b\pgfmathresult;%
\pgfmathparse{int(#1-2)}\let\a\pgfmathresult;%
\foreach \x in {1,...,\a} { \node (e\x) at (1.7*\x,0) {$\bullet$}; }%
\node (e\b) at (1.7*\a+1.2,1.2) {$\bullet$};
\node (e\c) at (1.7*\a+1.2,-1.2) {$\bullet$};

\foreach \x in {2,...,\a} {%
\pgfmathparse{int(\x-1)}\let\y\pgfmathresult;%
\draw[-] (e\x)--(e\y) node[mk@vertex, pos=.5,sloped,above] {\foreach
\i/\j/\valA/\valB in #2
{\ifthenelse{\IfVal{\j}{\x}{\i}{\y}}{$(\valA,\valB)$}{}}}; }
\draw[-] (e\c)--(e\a);
\draw[-] (e\b)--(e\a);
\end{tikzpicture}
}

\newcommand{\CanA}[2]
{
\begin{tikzpicture}[scale=.666]
\pgfmathparse{int(#1)}\let\x\pgfmathresult;%
\pgfmathparse{int(#1-1)}\let\y\pgfmathresult;%
\pgfmathparse{int(#1-2)}\let\z\pgfmathresult;%

\node (c1) at (\x,0) {$\bullet$};
\node (c\y) at (0,0) {$\bullet$};
\node (c\x) at (\x/2,-1.2) {$\bullet$};

\foreach \i in {2,...,\z} { \pgfmathparse{int(\i-1)}\let\a\pgfmathresult; \node (c\i) at (\x-\a*\x/\z,1.2) {$\bullet$}; }%

\foreach \i in {1,...,\z} {\pgfmathparse{int(\i+1)}\let\a\pgfmathresult; \draw[-latex] (c\i)--(c\a);}
\draw[-latex] (c1)--(c\x) node[mk@vertex, pos=.5,sloped,above] {\ifthenelse{#2=0}{$(a,b)$}{$(b,a)$}};
\draw[-latex] (c\x)--(c\y) node[mk@vertex, pos=.5,sloped,above] {\ifthenelse{#2=0}{$(b,a)$}{$(a,b)$}};

\end{tikzpicture}
} 

\newtheorem{thm}{Theorem}[section]
\newtheorem{lem}[thm]{Lemma}

\newtheorem{prop}[thm]{Proposition}
\newtheorem{rmk}[thm]{Remark}
\newtheorem*{pf}{Proof}

\newcommand{\nohyphens}[0]{\hyphenpenalty=10000\exhyphenpenalty=10000\relax}
\newcommand{\rhyphens}[0]{\hyphenpenalty=50\exhyphenpenalty=50\relax}

\begin{document}


\title{Self-injective artin algebras without short cycles in the component quiver}
\author{Maciej Karpicz}

\maketitle

\begin{abstract}
\nohyphens We give a complete description of all self-injective artin algebras of infinite representation type whose component quiver has no short cycles. \rhyphens
\end{abstract}

\footnote{e-mail: mkarpicz@mat.umk.pl}
\footnote{Faculty of Mathematics and Computer Science, Nicolaus Copernicus University, Chopina 12/18, 87-100 Toru{\'n}, Poland}


\section{Introduction and the main result}\label{sec1}
Throughout the paper, by an algebra we mean a~basic, connected, artin algebra over a~commutative artinian ring $k$.
For an algebra $A$, we denote by $\modd A$ the category of finitely generated right $A$-modules and by $\ind A$ the full subcategory of $\modd A$ given by the indecomposable modules.
An algebra $A$ is called self-injective if $A_A$ is an injective module, or equivalently, the projective and injective modules in $\modd A$ coincide.
A~prominent class of self-injective algebras  is formed by the orbit algebras $\hB/G$, where $\hB$ is the repetitive category of an algebra $B$ and $G$ is an admissible group of automorphisms of $\hB$.

An important combinatorial and homological invariant of the module category $\modd A$ of an algebra $A$ is its Auslander\dywiz Reiten quiver $\Gamma_A=\Gamma(\modd A)$.
It describes the structure of the quotient category $\modd A/\rad^\infty_A$, where $\rad^\infty_A$ is the infinite Jacobson radical of $\modd A$.
By a~result of Auslander \cite{Au}, $A$ is of finite representation type if and only if $\rad^\infty_A=0$.
Recall that $A$ is said to be of finite representation type if $\modd A$ admits only finitely many pairwise non-isomorphic indecomposable modules.

In general, it is important to study the behavior of the components of $\Gamma_A$ in the category $\modd A$.
Following \cite{S2}, a~component $\mc{C}$ of $\Gamma_A$ is called {\it generalized standard} if $\rad^\infty(X,Y)=0$ for all modules $X$ and $Y$ in $\mc{C}$.
It has been proved in \cite{S2} that every generalized standard component $\mc{C}$ of $\Gamma_A$ is {\it almost periodic}, that is, all but finitely many $\tau_A$-orbits in $\mc{C}$ are periodic.
Moreover, by a~result of \cite{SZ}, the additive closure $\add(\mc{C})$ of a~generalized standard component $\mc{C}$ of $\Gamma_A$ is closed under extensions in $\modd A$.
A~component of $\Gamma_A$ of the form $\bbZ\bbA_\infty/(\tau_A^r)$, where $r$ is a~positive integer, is called a~{\it stable tube} of rank $r$.
We note that, for $A$ self-injective, every infinite, generalized standard component $\mc{C}$ of $\Gamma_A$ is either acyclic with finitely many $\tau_A$-orbits or is a~{\it quasi-tube} (the stable part $\mc{C}^s$ of $\mc{C}$ is a~stable tube).
Following \cite{S9}, a~{\it component quiver} $\Sigma_A$ of an algebra $A$ has the components of $\Gamma_A$ as the vertices and two components $\mc{C}$ and $\mc{D}$ of $\Gamma_A$ are linked in $\Sigma_A$ by an arrow $\mc{C}\to\mc{D}$ if $\rad_A^\infty(X,Y)\neq 0$ for some modules $X$ in $\mc{C}$ and $Y$ in $\mc{D}$.
In particular, a~component $\mc{C}$ of $\Gamma_A$ is generalized standard if and only if $\Sigma_A$ has no loop at $\mc{C}$.
By a~{\it short cycle} in $\Sigma_A$ we mean a~cycle $\mc{C}\to\mc{D}\to\mc{C}$, where possibly $\mc{C}=\mc{D}$. We also mention that by a~result in \cite{JS} the component quiver $\Sigma_A$ of a~self-injective algebra $A$ of infinite representation type is fully cyclic, that is, any finite number of components of $\Gamma_A$ lies on a~common cycle in $\Sigma_A$.

For a~field $k$ of characteristic different from $2$, the {\it exceptional tubular algebra} $B_{ex}$ is the tubular algebra of the tubular type $(2,2,2,2)$, which is given by the following ordinary quiver
\begin{center}
\begin{tikzcd}
1 & 3\arrow{l}[swap]{\alpha}\arrow{dl}[near start, above=2]{\sigma} & 5\arrow{l}[swap]{\zeta}\arrow{dl}[near start, above=2]{\eta}\\
2 & 4\arrow{l}{\beta}\arrow{ul}[near start, below=2]{\gamma} & 6\arrow{l}{\omega}\arrow{ul}[near start, below=2]{\mu}
\end{tikzcd}
\end{center}
and the relations $\zeta\alpha=\eta\gamma$, $\mu\alpha=\omega\gamma$, $\zeta\sigma=\eta\beta$ and $\mu\sigma=-\omega\beta$.
Moreover, an automorphism $\varphi$ of the exceptional tubular algebra $B_{ex}$ is said to be {\it distinguished} if $\varphi(\gamma)=a\sigma$, $\varphi(\sigma)=b\gamma$, $\varphi(\beta)=c\alpha$, $\varphi(\alpha)=d\beta$, $\varphi(\mu)=e\eta$, $\varphi(\eta)=r\mu$, $\varphi(\omega)=u\zeta$ and $\varphi(\zeta)=v\omega$ for $a$, $b$, $c$, $d$, $e$, $r$, $u$, $v\in k\setminus\{0\}$ satisfying the following relations $dv=-ar$, $de=au$, $bv=cr$ and $be=-cu$.
In fact, $a=-\frac{dv}{r}$, $b=\frac{cr}{v}$, $e=-\frac{uv}{r}$, for $c$,~$d$,~$r$,~$u$,~$v\in k\setminus\{0\}$.

The aim of the paper is to prove the following theorem characterizing the class of representation-infinite self-injective artin algebras whose component quiver $\Sigma_A$ contains no short cycles.

\begin{thm}\label{thm_main}
Let $A$ be a~basic, connected, self-injective artin algebra of infinite representation type over an artinian ring $k$.
The following statements are equivalent.
\begin{enumerate}[(i)]
    \item The component quiver $\Sigma_A$ has no short cycles.
    \item $k$ is a~field and $A$ is isomorphic to an orbit algebra $\hB/G$, where $B$ is a~tilted algebra of Euclidean type or a~tubular algebra over $k$ and $G$ is an infinite cyclic group of automorphisms of $\hB$ of one of the following forms:
\begin{enumerate}[(a)]
\item $G=(\varphi\nu_{\hB}^2)$, for a~strictly positive automorphism $\varphi$ of $\hB$,
\item $G=(\varphi\nu_{\hB}^2)$, for an exceptional tubular algebra $B$ and a~rigid automorphism $\varphi$ of $\hB$, whose restriction to $B$ is a~distinguished automorphism of $B$,
\end{enumerate}
where $\nu_{\hB}$ is the Nakayama automorphism of $\hB$.
\end{enumerate}
\end{thm}

The assumptions of the theorem imply that the module category $\modd A$ of $A$ contains no infinite short cycles and every component $\mc{C}$ in $\Gamma_A$ has no external short paths.
By a~{\it short cycle} in $\modd A$ we mean a~sequence $M\stackrel{f}{\to}N\stackrel{g}{\to}M$ of non-zero non-isomorphisms between indecomposable modules in $\modd A$ \cite{RSS}, and such a~cycle is said to be {\it infinite} if at least one of the homomorphisms $f$ or $g$ belongs to $\rad^\infty_A$.
Moreover, following \cite{RS}, by an {\it external short path} of a~component $\mc{C}$ of $\Gamma_A$ we mean a~sequence $X\to Y\to Z$ of non-zero homomorphisms between indecomposable modules in $\modd A$ with $X$ and $Z$ in $\mc{C}$ but $Y$ not in $\mc{C}$.

The paper is organized as follows.
In Section $2$ we recall the essential background.
In Section $3$ we show that every automorphism of non-exceptional tubular algebra fixes a~point.
Finally, Section $4$ is devoted to the proof of the Theorem \ref{thm_main}.

For a~basic background on the representation theory of algebras applied in the paper we refer to the books \cite{ASS}, \cite{R2}, \cite{SS1}, \cite{SS2}, \cite{SY} and to the survey articles \cite{MS3}, \cite{S7} and \cite{SY6}. 

\section{Preliminaries}\label{sec2}
Let $A$ be an artin algebra over a commutative artin ring $k$.
We denote by $\Gamma_A$ the Auslander\dywiz Reiten quiver of $A$.
Recall that $\Gamma_A$ is a~valued translation quiver whose vertices are the isomorphism classes $\{X\}$ of modules $X$ in $\ind A$, the valued arrows of $\Gamma_A$ describe minimal left almost split morphisms with indecomposable domain and minimal right almost split morphisms with indecomposable codomain, and the translation is given by the Auslander\dywiz Reiten translations $\tau_A=\dtr$ and $\tau_A^{-1}=\trd$ (see \cite[Chapter III]{SY} for details).
We will identify vertices of $\Gamma_A$ with the corresponding indecomposable modules and by a~component in $\Gamma_A$ we mean a~connected component.
Let $\mc{C}$ be a~family of components of $\Gamma_A$.
Then $\mc{C}$ is said to be {\it sincere} if any simple $A$-module occurs as a~composition factor of a~module in $\mc{C}$, and {\it faithful} if its annihilator $\ann_A(\mc{C})$ in $A$ (the intersection of the annihilators of all modules in $\mc{C}$) is zero.
Observe that if $\mc{C}$ is faithful, then $\mc{C}$ is sincere.
Moreover, the family $\mc{C}$ is said to be {\it separating} in $\modd A$ if the indecomposable modules in $\modd A$ split into three disjoint classes $\mc{P}^A$, $\mc{C}^A=\mc{C}$ and $\mc{Q}^A$ such that:
\begin{enumerate}
    \item[(S1)] $\mc{C}^A$ is a~sincere generalized standard family of components;
    \item[(S2)] $\Hom_A(\mc{Q}^A,\mc{P}^A)=0$, $\Hom_A(\mc{Q}^A,\mc{C}^A)=0$, $\Hom_A(\mc{C}^A,\mc{P}^A)=0$;
    \item[(S3)] any homomorphism from $\mc{P}^A$ to $\mc{Q}^A$ factors through the additive category $\add\mc{C}^A$ of $\mc{C}^A$.
\end{enumerate}
Moreover, a~separating family $\mc{C}^A=(\mc{C}^A_i)_{i\in I}$ is {\it strongly separating} if
\begin{enumerate}
    \item[(S4)] any homomorphism from $\mc{P}^A$ to $\mc{Q}^A$ factors trough the additive category $\add\mc{C}^A_i$, for any $i\in I$.
\end{enumerate}
Let $\Lambda$ be a~canonical algebra in the sense of Ringel \cite{R2} (and \cite{R3}).
Then the quiver $Q_\Lambda$ of $\Lambda$ has a~unique sink and a~unique source.
Denote by $Q^*_\Lambda$ the quiver obtained from $Q_\Lambda$ by removing the unique source of $Q_\Lambda$ and the arrows attached to it.
Then $\Lambda$ is said to be a~canonical algebra of Euclidean type (respectively, of tubular type) if $Q^*_\Lambda$ is a~ Dynkin quiver (respectively, a~Euclidean quiver).
The general shape of the Auslander-Reiten quiver $\Gamma_\Lambda$ of $\Lambda$ is as follows: $$\Gamma_\Lambda=\mc{P}^\Lambda\vee\mc{T}^\Lambda\vee\mc{Q}^\Lambda,$$ where $\mc{P}^\Lambda$ is a~family of components containing a~unique preprojective component $\mc{P}(\Lambda)$ and all indecomposable projective $\Lambda$-modules, $\mc{Q}^\Lambda$ is a~family of components containing a~unique preinjective component $\mc{Q}(\Lambda)$ and all indecomposable injective $\Lambda$-modules, and $\mc{T}^\Lambda$ is an infinite family of pairwise orthogonal, generalized standard, faithful stable tubes, separating $\mc{P}^\Lambda$ from $\mc{Q}^\Lambda$, and with all but finitely many stable tubes of rank one.
An algebra $C$ of the form $\End_\Lambda(T)$, where $T$ is a~multiplicity-free tilting module from the additive category $\add(\mc{P}^\Lambda)$ of $\mc{P}^\Lambda$ is said to be a~{\it concealed canonical algebra} of type $\Lambda$ (\cite{LP1}).
More generally, an algebra $B$ of the form $\End_\Lambda(T)$, where $T$ is a~multiplicity-free tilting module from the additive  category $\add(\mc{P}^\Lambda\cup\mc{T}^\Lambda)$ of $\mc{P}^\Lambda\cup\mc{T}^\Lambda$ is said to be an almost concealed canonical algebra of type $\Lambda$ \cite{LS1}.
An almost concealed canonical algebra $B$ of a~tubular type is called a~{\it tubular algebra}.
Recall that a~tubular algebra $B$ may be obtained as a~tubular branch extension of a~concealed algebra of Euclidean type(see \cite{R2} and \cite{R3}).
Moreover, an almost concealed canonical algebra of Euclidean type is a~representation-infinite, tilted algebra of Euclidean type for which the preinjective component contains all indecomposable injective modules (see \cite{R2}).

For an algebra $A$, we denote by $\D$ the standard duality $\Hom_k(-,E)$ on $\modd A$, where $E$ is a~minimal injective cogenerator in $\modd k$.
Then an algebra $A$ is {\it self-injective} if and only if $A\cong\D(A)$ in $\modd A$.

Let $A$ be a~self-injective algebra and $1_A=\{e_i\;|\;1\leq i\leq s\}$ a~complete set of pairwise orthogonal primitive idempotents of $A$ whose sum is the identity $1_A$ of $A$.
We denote by $\nu=\nu_A$ the {\it Nakayama automorphism} of $A$ inducing an $A-A$-bimodule isomorphism $A\cong \D(A)_\nu$, where $\D(A)_\nu$ denotes the right $A$-module obtained from $\D(A)$ by changing the right operation of $A$ as follows: $f\cdot a=f\nu(a)$ for each $a\in A$ and $f\in\D(A)$.
Hence we have $\soc(\nu(e_i)A)\cong\top(e_iA)$ $(=e_i A/\rad(e_i A))$ as right $A$-modules for all $i\in\{1,\ldots,s\}$.
Since $\{\nu(e_i)A\;|\; 1\leq i\leq s\}$ is a~complete set of representatives of indecomposable projective right $A$-modules, there is a~{\it (Nakayama) permutation} of $\{1,\ldots,s\}$, denoted again by $\nu$, such that $\nu(e_i)A\cong e_{\nu(i)}A$ for all $i\in\{1,\ldots,s\}$.
Invoking the Krull-Schmidt theorem, we may assume that $\nu(e_i A)=\nu(e_i)A=e_{\nu(i)}A$ for all $i\in\{1,\ldots,s\}$.

Let $B$ be an algebra and $1_B=e_1+\cdots+e_n$ a~decomposition of the identity of $B$ into the sum of pairwise orthogonal primitive idempotents of $B$.
We associate to $B$ a~self-injective, locally bounded $k$-category $\hB$, called the {\it repetitive category} of $B$ (see \cite{HW}). The objects of $\hB$ are $e_{m,i}$, $m\in\bbZ$, $i\in\{1,\ldots,n\}$ and the morphism spaces are defined in the following way
$$
\hB(e_{m,i},e_{r,j})=\left\{\begin{matrix}
e_jBe_i,&r=m,\\
D(e_iBe_j),&r=m+1,\\
0,&\mbox{otherwise}.
\end{matrix}\right.
$$
Observe that $e_jBe_i=\Hom_B(e_iB,e_jB)$, $D(e_iBe_j)=e_jD(B)e_i$ and
$$\bigoplus_{(r,i)\in\bbZ\times\{1,\ldots,n\}}\hB(e_{m,i},e_{r,j})=e_jB\oplus\D(Be_j),$$
for any $m\in\bbZ$ and $j\in\{1,\ldots,n\}$.
We denote by $\nu_{\hB}$ the {\it Nakayama automorphism} of $\hB$ defined by $$\nu_{\hB}(e_{m,i})=e_{m+1,i},\mbox{ for any }(m,i)\in\bbZ\times\{1,\ldots,n\}.$$
An automorphism $\varphi$ of the $k$-category $\hB$ is said to be:
\begin{itemize}
    \item {\it positive} if for each pair $(m,i)\in\bbZ\times\{1,\ldots,n\}$ we have $\varphi(e_{m,i})=e_{p,j}$ for some $p\geq m$ and $j\in\{1,\ldots,n\}$;
    \item {\it rigid} if for each pair $(m,i)\in\bbZ\times\{1,\ldots,n\}$ there exists $j\in\{1,\ldots,n\}$ such that $\varphi(e_{m,i})=e_{m,j}$;
    \item {\it strictly positive} if $\varphi$ is positive but not rigid.
\end{itemize}

A~group $G$ of automorphisms of $\hB$ is said to be {\it admissible} if $G$ acts freely on the set of objects of $\hB$ and has finitely many orbits.
Then we may consider the orbit category $\hB/G$ of $\hB$ with respect to $G$ whose objects are $G$-orbits of objects in $\hB$, and the morphism spaces are given by
$$
(\hB/G)(a,b)=\left\{f_{y,x}\in\prod_{(x,y)\in a\times b}\hB(x,y)\mid \operatorname{\forall}_{g\in G,(x,y)\in a\times b}\;gf_{y,x}=f_{gy,gx}\right\}
$$
for all objects $a$,~$b$ of $\hB/G$.
Since $\hB/G$ has finitely many objects and the morphism spaces in $\hB/G$ are finitely generated, we have the associated self-injective artin algebra $\bigoplus(\hB/G)$ which is the direct sum of all morphism spaces in $\hB/G$, called the {\it orbit algebra} of $\hB$ with respect to $G$.
For example, the infinite cyclic group $(\nu_{\hB})$ generated by $\nu_{\hB}$ is admissible and $\hB/(\nu_{\hB})$ is the trivial extension $B\ltimes\D(B)$ of $B$ by $\D(B)$.

We refer to \cite{SY2} and \cite{SY5} for criteria on a~self-injective algebra $A$ to be an orbit algebra $\hB/G$ with $G$ an infinite cyclic group generated by an automorphism of the form $(\varphi\nu_{\hB})$, with $\varphi$ a~strictly positive automorphism of $\hB$. 

\section{Tubular algebras}\label{sec4}
From now on we assume that an automorphism of an algebra $B$ is rigid.
Let $1_B=e_1+\cdots+e_n$ be a~decomposition of the identity of $B$ into the sum of pairwise orthogonal, primitive idempotents.
We say that an automorphism $\varphi$ of $B$ is {\it rigid} if $\varphi(\{e_1,\ldots,e_n\})=\{e_1,\ldots,e_n\}$.

The aim of this section is to show that for a~tubular algebra $B$ and a rigid automorphism $\varphi$ of $B$, if $B$ is not an exceptional algebra and $\varphi$ a~distinguished automorphism, then $\varphi$ admits a~fixed point, that is for the induced, by $\varphi$, automorphism $\varphi$ of the ordinary quiver $Q_B$ of the algebra $B$, there is a~vertex $e$ such that $\varphi(e)=e$.
This will be a~consequence of Lemma \ref{lm_3.1} and Proposition \ref{prop_3.2}


Let $B$ be a~tubular algebra which is a~tubular branch extension of a~concealed algebra of Euclidean type $\wt{A}_n$. If $\varphi$ is an automorphism of the ordinary quiver $Q_B$ of $B$, then $\varphi$ acts freely on the set of vertices of $Q_B$ only if $Q_B$ is the ordinary quiver of the exceptional tubular algebra $B_{ex}$ (see \cite[Section 3]{S1}).
Therefore, we start with the following lemma describing the case of an exceptional tubular algebra $B_{ex}$.
\begin{lem}\label{lm_3.1}
Let $B$ be a~tubular algebra over a~field $k$ given by the quiver
\begin{center}
\begin{tikzcd}
1 & 3\arrow{l}[swap]{\alpha}\arrow{dl}[near start, above=2]{\sigma} & 5\arrow{l}[swap]{\zeta}\arrow{dl}[near start, above=2]{\eta}\\
2 & 4\arrow{l}{\beta}\arrow{ul}[near start, below=2]{\gamma} & 6\arrow{l}{\omega}\arrow{ul}[near start, below=2]{\mu}
\end{tikzcd}
\end{center}
with relations $\alpha\zeta=\gamma\nu$, $\alpha\mu=\gamma\omega$, $\sigma\zeta=\beta\nu$ and $\sigma\mu=x\beta\omega$, where $x\in k\setminus\{0,1\}$. Let $\varphi$ be an automorphism of $B$. Then $\varphi$ acts freely on vertices of $Q_B$ if and only if $B$ is exceptional and $\varphi$ a~distinguished automorphism of $B$.
\end{lem}
\begin{pf}
Taking the values of the automorphism $\varphi$ of $B$ on the above equalities we get $dv\beta\omega=ar\sigma\mu$, $de\beta\eta=au\sigma\zeta$, $bv\gamma\omega=cr\alpha\mu$, $be\gamma\eta=xcu\alpha\zeta$, and hence the equalities $dv=xar$, $de=au$, $bv=cr$ and $be=xcu$.
Therefore, combining those equalities, we get $xaer=dev=auv$ and $cer=bev=xcuv$, which implies $xer=uv$, $er=xuv$, hence $x^2=1$.
Because $x\neq 0,1$, we get a~contradiction if and only if $x\neq -1$ and $k$ is not a~field of characteristic $2$.
Therefore, $\varphi$ does not fix a~point if $x=-1\neq 1$, that is $B$ is the exceptional tubular algebra and $\varphi$ a~distinguished automorphism of $B$.\hfill{$\square$}
\end{pf}

\begin{rmk}
The above lemma corrects \cite[Lemma 3.5]{S1}.
\end{rmk}

By {\it small Euclidean graphs} we understand the following Euclidean graphs: $\wt{A}_{11}$, $\wt{A}_{12}$, $\wt{B}_m$, $\wt{C}_m$, $\wt{BC}_m$, for $m=2,\ldots,5$, $\wt{BD}_m$, $\wt{CD}_m$, for $m=3,4,5$, $\wt{F}_{41}$, $\wt{F}_{42}$, $\wt{G}_{21}$ and $\wt{G}_{22}$ (see \cite{DR}).

Let $H$ be a~representation-infinite hereditary artin algebra.
Recall that, following \cite[Section XIV.4]{SS1}, a~{\it concealed domain} $\mc{DP}(H)$ of $H$ is a~subquiver of the postprojective component $\mc{P}(H)$ of $\Gamma_H$ with the following properties:
\begin{enumerate}[(d1)]
\item $\mc{DP}(H)$ is a~finite full translation subquiver of
$\mc{P}(H)$ which is closed under the predecessors in
$\mc{P}(H)$,
\item for any multiplicity-free postprojective tilting $H$-module $T=T_1\oplus\cdots\oplus T_n$, there exists a~ postprojective tilting module $T'=T'_1\oplus\cdots\oplus T'_n$ such that the $H$-modules $T_1',\ldots,T'_n$ are indecomposable, pairwise non-isomorphic, lie in $\mc{DP}(H)$, and there is an isomorphism of algebras $$\End T_H\cong\End T'_H.$$
\end{enumerate}
We set $\mc{DP}(H):=\{\tau^{-r}_HP(a)\mid r\in\{1,\ldots,r_\Delta\}\mbox{ and }a\in\Delta_0\}$, where $r_\Delta$ is a~ minimal positive integer such that $\Hom_H(P(a),\tau^{-r}_HP(b))\neq 0$ for all $r\geq r_\Delta$ and $a,b\in\Delta_0$.
Note, that an integer $r_\Delta$ exists because the postprojective component $\mc{P}(H)$ of a~hereditary algebra $H$, contains only finitely many non-sincere modules.
Therefore, in order to compute the Gabriel quivers of concealed algebras, it is sufficient to calculate the Gabriel quivers of a~tilted algebras $\End_H(T)$, for which the tilting module $T$ has all direct summands from $\mc{DP}(H)$.

\begin{rmk}\label{rmk_1}
Note that for a~given hereditary algebra $H$ of Euclidean type, different from $\wt{A}_n$ type, the concealed domain $\mc{DP}(H)$ of $H$ admits a~tilting module $T$ such that the quiver of $B=\End_H(T)$ has the same underlying graph as $\Delta$ and an arbitrary orientation of arrows.
Moreover, there is an equivalence of categories $\mc{F}(T):=\{Y\in\modd B\mid\Ext_B^1(U,\D(T))=0\}$ in $\modd H$, containing $\mc{DP}(H)$, and $\mc{X}(T)=\{X\in\modd B\mid\Hom_B(X,\D(T))=0\}$ in $\modd B$.
\end{rmk}

\begin{prop}\label{prop_3.1}\label{prop_3.2}
Let $B$ be a~non-exceptional tubular algebra. Then any automorphism $\varphi\in\End_k(B)$ fixes a~point of $Q_B$.
\end{prop}
\begin{pf}
Let $C$ be a~tame concealed algebra of Euclidean type $\Delta$ such that $B$ is a~branch extension of $C$ of tubular type.

Clearly, for any automorphism $\varphi\colon B\to B$ its restriction to $C$ is an automorphism of $C$.
Therefore, to prove that $\varphi$ fixes a~point in $B$, it suffices to show that the restriction $\varphi\mid_C$ of $\varphi$ to $C$ fixes a~point.
In order to show this, we need shapes of ordinary quivers of tame concealed algebras.

The classification of concealed algebras of Euclidean type $\Delta$, where $\Delta$ is one of Euclidean graphs $\wt{A}_n$, for $n\geq 2$, $\wt{D}_n$, for $n\geq 4$, $\wt{E}_6, \wt{E}_7, \wt{E}_8$, in terms of quivers with relations was given by Bongartz \cite{Bo2} and Happel-Vossieck \cite{HV} (see also \cite[Section XIV]{SS1}).
Moreover, a~simple inspection of the Bongartz-Happel-Vossieck list shows that every automorphism of a~concealed algebra of Euclidean type, different from $\wt{A}_n$ type, fixes a~point.

Now, because tubular algebras also arise from tilting modules over canonical algebras of tubular type, then it follows that the Grothendieck group has a~small rank (see \cite[Theorem 3.2]{SY6}).
Therefore, we will provide possible shapes of ordinary quivers of concealed algebras of Euclidean type of the form $C=\End_H(T)$, where $T$ is a~multiplicity-free tilting module and $H$ is a~hereditary algebra, whose ordinary quiver is one of oriented small Euclidean graphs (see \cite{DR} for realization of such hereditary algebras).

By Remark \ref{rmk_1} we may consider only hereditary algebras whose ordinary quiver, denoted by $\Delta(G)$, where $G$ is one of small Euclidean graphs, is equipped with an orientation of arrows from the right hand side to the left hand side.
Computation of tilting modules in the concealed domains of such hereditary algebras is rather an easy, but tedious, task and was done by a~computer (see \cite{K4}).
We list here only the shapes (frames) of ordinary quivers of these concealed algebras.
In the list of frames below, an unoriented arrow may have an arbitrary orientation.

\begin{enumerate}
\item For a~hereditary algebra $H$ whose ordinary quiver is $\Delta(\wt{B}_2)$ (respectively, $\Delta(\wt{A}_{11})$, $\Delta(\wt{A}_{12})$, $\Delta(\wt{B}_3)$, $\Delta(\wt{B}_4)$, $\Delta(\wt{BC}_2)$, $\Delta(\wt{BC}_3)$, $\Delta(\wt{BC}_4)$, $\Delta(\wt{C}_2)$, $\Delta(\wt{C}_3)$ and $\Delta(\wt{C}_4)$) we get the frame $\wt{B}_2$ (respectively, $\wt{A}_{11}$, $\wt{A}_{12}$, $\wt{B}_3$, $\wt{B}_4$, $\wt{BC}_2$, $\wt{BC}_3$, $\wt{BC}_4$, $\wt{C}_2$, $\wt{C}_3$ and $\wt{C}_4$).

\item For a~hereditary algebra $H$ whose ordinary quiver is $\Delta(\wt{BD}_3)$ (respectively, $\Delta(\wt{CD}_3)$) we have the following frames:
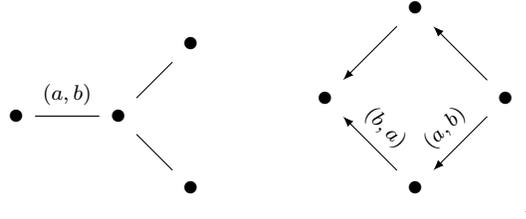
\begin{center}
\DGraph{4}{{1/2/a/b}}\hspace{1cm}
\begin{tikzpicture}
\node (c1) at (0,0) {$\bullet$};
\node (c2) at (1.2,1.2) {$\bullet$};
\node (c3) at (1.2,-1.2) {$\bullet$};
\node (c4) at (2.4,0) {$\bullet$};
\draw[-latex] (c2)--(c1);
\draw[-latex] (c4)--(c2);
\draw[-latex] (c3)--(c1) node[mk@vertex, pos=.5,sloped,above] {$(b,a)$};
\draw[-latex] (c4)--(c3) node[mk@vertex, pos=.5,sloped,above] {$(a,b)$};
\end{tikzpicture},
\end{center}
where $(a,b)=(1,2)$ (respectively, $(a,b)=(2,1)$).
\item For a~hereditary algebra $H$ whose ordinary quiver is $\Delta=\wt{BD}_4$ (respectively, $\wt{CD}_4$) we have the following frames:
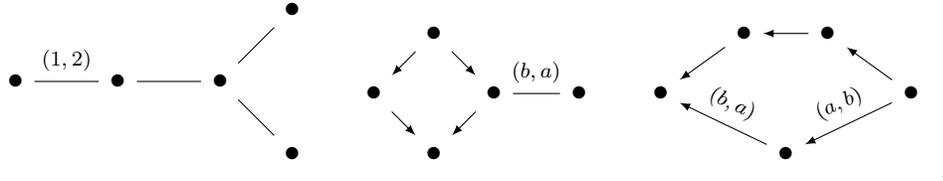
\begin{center}
\DGraph{5}{{1/2/1/2}}\hspace{.3cm}
\begin{tikzpicture}[scale=.666]
\node (a) at (0,0) {$\bullet$};
\node (b) at (1.2,1.2) {$\bullet$};
\node (c) at (1.2,-1.2) {$\bullet$};
\node (d) at (2.4,0) {$\bullet$};
\node (e) at (4.1,0) {$\bullet$};
\draw[-latex] (b)--(a);
\draw[-latex] (a)--(c);
\draw[-latex] (b)--(d);
\draw[-latex] (d)--(c);
\draw[-] (e)--(d) node[mk@vertex, pos=.5,sloped,above] {$(b,a)$};
\end{tikzpicture}\hspace{.3cm}
\CanA{5}{0},
\end{center}
where $(a,b)=(1,2)$ (respectively, $(a,b)=(2,1)$).

\item For a~hereditary algebra $H$ whose ordinary quiver is $\Delta(\wt{F}_{41})$ (respectively, $\Delta(\wt{F}_{42})$) we have the following frames:
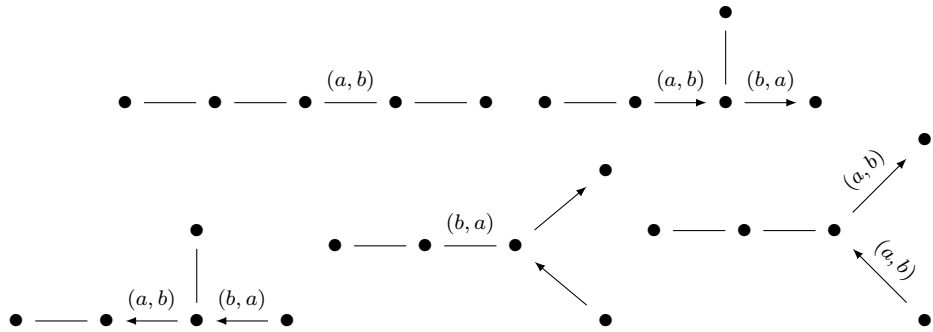
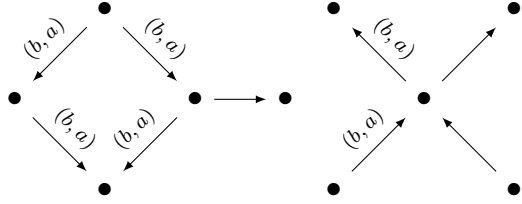
\begin{center}
\LinGraph{5}{{3/4/a/b}}
\begin{tikzpicture}
\foreach \x in {1,...,4} {\node (e\x) at (1.2*\x-1.2,0) {$\bullet$};}
\node (e5) at (2.4,1.2) {$\bullet$}; \draw[-] (e2)--(e1);
\draw[-latex] (e2)--(e3) node[mk@vertex, pos=.5,sloped,above] {$(a,b)$};
\draw[-latex] (e3)--(e4) node[mk@vertex, pos=.5,sloped,above] {$(b,a)$};
\draw[-] (e5)--(e3);
\end{tikzpicture}
\begin{tikzpicture}
\foreach \x in {1,...,4} {\node (e\x) at (1.2*\x-1.2,0) {$\bullet$};}
\node (e5) at (2.4,1.2) {$\bullet$}; \draw[-] (e1)--(e2);
\draw[-latex] (e3)--(e2) node[mk@vertex, pos=.5,sloped,above] {$(a,b)$};
\draw[-latex] (e4)--(e3) node[mk@vertex, pos=.5,sloped,above] {$(b,a)$};
\draw[-] (e5)--(e3);
\end{tikzpicture}
\begin{tikzpicture}
\node (e1) at (0,0) {$\bullet$};
\node (e2) at (1.2,0) {$\bullet$};
\node (e3) at (2.4,0) {$\bullet$};
\node (e4) at (3.6,1) {$\bullet$};
\node (e5) at (3.6,-1) {$\bullet$};
\draw[-] (e1)--(e2);
\draw[-] (e2)--(e3) node[mk@vertex, pos=.5,sloped,above] {$(b,a)$};
\draw[-latex] (e3)--(e4);
\draw[-latex] (e5)--(e3);
\end{tikzpicture}
\begin{tikzpicture}
\foreach \x in {1,...,3} {\node (e\x) at (1.2*\x-1.2,0) {$\bullet$};}
\node (e4) at (3.6,1.2) {$\bullet$};
\node (e5) at (3.6,-1.2) {$\bullet$};
\draw[-] (e1)--(e2);
\draw[-] (e2)--(e3);
\draw[-latex] (e3)--(e4) node[mk@vertex, pos=.5,sloped,above] {$(a,b)$};
\draw[-latex] (e5)--(e3) node[mk@vertex, pos=.5,sloped,above] {$(a,b)$};
\end{tikzpicture}
\end{center}
\begin{center}
\begin{tikzpicture}
\node (e1) at (1.2,1.2) {$\bullet$}; \node (e2) at (0,0) {$\bullet$};
\node (e3) at (2.4,0) {$\bullet$}; \node (e4) at (1.2,-1.2) {$\bullet$};
\node (e5) at (3.6,0) {$\bullet$};
\draw[-latex] (e1)--(e2) node [mk@vertex, pos=.5,sloped,above] {$(b,a)$};
\draw[-latex] (e2)--(e4) node [mk@vertex, pos=.5,sloped,above] {$(b,a)$};
\draw[-latex] (e1)--(e3) node [mk@vertex, pos=.5,sloped,above] {$(b,a)$};
\draw[-latex] (e3)--(e4) node [mk@vertex, pos=.5,sloped,above] {$(b,a)$};
\draw[-latex] (e3)--(e5);
\end{tikzpicture}\vspace{1cm}
\begin{tikzpicture}
\node (e1) at (0,1.2) {$\bullet$};
\node (e2) at (0,-1.2) {$\bullet$};
\node (e3) at (1.2,0) {$\bullet$};
\node (e4) at (2.4,1.2) {$\bullet$};
\node (e5) at (2.4,-1.2) {$\bullet$};
\draw[-latex] (e3)--(e1) node[mk@vertex, pos=.5,sloped,above] {$(b,a)$};
\draw[-latex] (e3)--(e4);
\draw[-latex] (e2)--(e3) node[mk@vertex, pos=.5,sloped,above] {$(b,a)$};
\draw[-latex] (e5)--(e3);
\end{tikzpicture}
\end{center}
where $(a,b)=(1,2)$ (respectively, $(a,b)=(2,1)$).
\item For a~hereditary algebra $H$ whose ordinary quiver is $\Delta(\wt{G}_{21})$ (respectively, $\Delta(\wt{G}_{22})$) we have the following two frames:
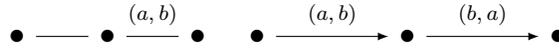
\begin{center}
\LinGraph{3}{{2/3/a/b}}
\begin{tikzpicture}
\foreach \x in {0,1,2} {\node (e\x) at (2*\x,0) {$\bullet$};}
\draw[-latex] (e0)--(e1) node[mk@vertex, pos=.5,sloped,above] {$(a,b)$};
\draw[-latex] (e1)--(e2) node[mk@vertex, pos=.5,sloped,above] {$(b,a)$};
\end{tikzpicture}
\end{center}
where $(a,b)=(1,3)$ (respectively, $(a,b)=(3,1)$).
\end{enumerate}
Now, a~simple inspection of the above listed frames shows that every automorphism of a~concealed algebra of Euclidean type $\Delta$, where $\Delta$ is one of small Euclidean graphs, has a~fixed point.\hfill{$\square$}
\end{pf} 

\section{Proof of the Theorem}\label{sec5}
Let $A$ be a self-injective artin algebra of infinite representation type such that the component quiver $\Sigma_A$ of $A$ contains no short cycles.
Then the Auslander\dywiz Reiten quiver $\Gamma_A$ of $A$ consists of modules which do not lie on infinite short cycles and all components in $\Gamma_A$ are generalized standard.
In particular, for any indecomposable $A$-module $M$, we have $\rad_A^\infty(M,M)=0$.

Given a~module $M$ in $\modd A$, we denote by $[M]$ the image of $M$ in the Grothendieck group $K_0(A)$ of $A$.
Thus $[M]=[N]$ if and only if the modules $M$ and $N$ have the same composition factors including the multiplicities.
We also mention that, by a~result proved in \cite{RSS}, every indecomposable module $M$ in $\modd A$ which does not lie on a~short cycle is uniquely determined by $[M]$ (up to isomorphism).
In addition, recall that, following \cite{S4}, a~family $\mc{C}=(\mc{C}_i)_{i\in I}$ of components of $\Gamma_A$ is said to have {\it common composition factors}, if, for each pair $i$ and $j$ in $I$, there are modules $X_i\in\mc{C}_i$ and $X_j\in\mc{C}_j$ with $[X_i]=[X_j]$.
Moreover, $\mc{C}$ is {\it closed under composition factors} if, for every indecomposable modules $M$ and $N$ in $\modd A$ with $[M]=[N]$, $M\in\mc{C}$ forces $N\in\mc{C}$.

We start with proving the following proposition which, for an algebraically closed field, was proved by Ringel in \cite[(5.2)]{R2}.
\begin{prop}\label{prop_4.1}
Let $B$ be a~tubular algebra with the canonical decomposition $$\Gamma_B=\mc{P}\vee\mc{T}_0\vee\bigvee_{q\in\bbQ^+}\mc{T}_q\vee\mc{T}_\infty\vee\mc{Q}$$ of the Auslander\dywiz Reiten quiver.
Then, for any $q\in\bbQ^{+}\cup\{0,\infty\}$, the family $\mc{T}_q$ of tubes is closed under composition factors.
\end{prop}
\begin{pf}
Let $M$ be a~$B$-module from $\mc{T}^B_p$ and $N$ be a~$B$-module from $\mc{T}^B_q$, $p$,~$q\in\bbQ^+\cup\{0,\infty\}$.
Assume that $[M]=[N]$.
We will show that $p=q$.
Assume to the contrary that $p<q$.
Take some $s\in\bbQ^+$ with $p<s<q$.
Since the family of stable tubes $\mc{T}_s^B=(\mc{T}_{s,x}^B)_{x\in\bbX_s}$ is infinite, there is $x\in\bbX_s$ such that $\mc{T}_{s,x}^B$ is a~stable tube of rank one.
Take the unique mouth module $X$ in $\mc{T}_{s,x}^B$.
Clearly, $X=\tau_BX$.
We know that the family $\mc{T}_s^B$ strongly separates $\mc{X}=\mc{P}^B\vee\bigvee_{l<s}\mc{T}_l^B$ from $\mc{Y}=\bigvee_{l>s}\mc{T}_l^B\vee\mc{Q}^B$, that is, every homomorphism $f$ from $\add\mc{X}$ to $\add\mc{Y}$ factors through $\add\mc{T}_{s,y}^B$, for every $y\in\bbX_s$.
Observe that the injective hull $E_B(M)$ of $M$ is in $\add(\mc{T}_\infty\vee\mc{Q})$ and the projective cover $P_B(N)$ of $N$ is in $\add(\mc{P}\vee\mc{T}_0)$.
Therefore, $\Hom_B(M,\mc{T}_{s,x}^B)\neq0$ and $\Hom_B(\mc{T}_{s,x}^B,N)\neq 0$.
Hence, applying \cite[Lemma 3.9]{S3} $\Hom_B(M,X)\neq 0$ and $\Hom_B(X,N)\neq 0$, because $\mc{C}_{s,x}^B$ is of rank one.
Next, since $\mc{T}_s^B$ separates $\mc{X}$ from $\mc{Y}$, we have $\Hom_B(X,M)=0$ and $\Hom_B(N,X)=0$.
Further, since $[M]=[N]$, applying \cite[Proposition 4.1]{S3} (see also \cite{AR}) we have the equalities
$$
\begin{array}{rcl}
\abs{\Hom_B(X,M)}-\abs{\Hom_B(M,X)}&=&\abs{\Hom_B(X,M)}-\abs{\Hom_B(M,\tau_B X)}\\
&=&\abs{\Hom_B(X,N)}-\abs{\Hom_B(N,\tau_B X)}\\
&=&\abs{\Hom_B(X,N)}-\abs{\Hom_B(N,X)}
\end{array}
$$
where $\abs{V}$ is the length of a~$k$-module $V$.
Then $\Hom_B(X,M)=0$ and $\Hom_B(N,X)=0$ lead to a~desired contradiction
$$0>-|\Hom_B(M,X)|=|\Hom_B(X,N)|>0.$$
\phantom{.}\hfill$\square$
\end{pf}

\begin{prop}\label{lem_4.2}
Let $A$ be a~basic, connected self-injective algebra of infinite representation type such that the component quiver $\Sigma_A$ of $A$ contains no short cycles. Then the Auslander\dywiz Reiten quiver $\Gamma_A$ of $A$ admits a~family $\mc{C}=(\mc{C})_{x\in\bbX}$ of quasi-tubes having common composition factors, closed under composition factors and consisting of modules which do not lie on infinite short cycles in $\modd A$.
\end{prop}
\begin{pf}
Because $\Sigma_A$ contains no short cycles, then every component in $\Gamma_A$ is generalized standard.
Therefore, since $A$ is of infinite representation type, applying \cite[Corollary 4.4]{S8}, we conclude that there exists an ideal $I$ in $A$ such that $A'=A/I$ is tame concealed.
Thus, there exists a~family of stable tubes $\mc{T}^{A'}=(\mc{T}_x)_{x\in\bbX}$ in $\Gamma_{A'}$ with common composition factors.
In addition, (see \cite[Section 2]{KSY}), there is an infinite family $\mc{C}=(\mc{C}_x)_{x\in\bbX}$ of quasi-tubes in $\Gamma_A$ such that, for any $x\in\bbX$, $\mc{T}_x\subseteq\mc{C}_x$ and the equality holds for almost all $x\in\bbX$.
Obviously, because $\mc{T}^{A'}$ is a~family of stable tubes in $\Gamma_{A'}$ with common composition factors, the family $\mc{C}$ is a~family of quasi-tubes with common composition factors.
We claim that $\mc{C}$ is closed under composition factors.

Let $N$ be a~module in $\Gamma_A$ and $M$ a~module in $\mc{C}=(\mc{C}_x)_{x\in\bbX}$.
Assume that $[M]=[N]$.
We will show that $N$ belongs to $\mc{C}$.
Let $\mc{C}_y$, for some $y\in\bbX$, be the quasi-tube in the family $\mc{C}$ containing $M$.
Let $1_A=e+f$ be a~decomposition of $1_A$ into a~sum of two idempotents such that all direct summands of $eA/\rad(eA)$ are isomorphic to the composition factors of modules in $\mc{C}$, but the module $fA/\rad(fA)$ has no such direct summands.
Consider the quotient algebra $B=A/AfA$.
Then $\mc{C}_y$ is a~component in $\Gamma_B$.
Moreover, the $A$-module $N$ is also a~module over $B$.

Since $\mc{C}_y$ is a~generalized standard quasi-tube without external short paths, applying \cite[Theorem A]{MS2}, we conclude that $B$ is a~quasi-tube enlargement of a~concealed canonical algebra $C$ and there is a~separating family $\mc{C}^B$ of quasi-tubes in $\Gamma_B$ containing the quasi-tube $\mc{C}_y$.
In particular, we have a~decomposition $\Gamma_{B}=\mc{P}^{B}\vee\mc{C}^{B}\vee\mc{Q}^{B}$ with $\mc{C}^B$ separating $\mc{P}^B$ from $\mc{Q}^B$.

Therefore, by dual arguments, we may assume that $N$ belongs to $\mc{P}^B$.

From \cite[Theorem C]{MS2} there is a~unique maximal tubular coextension $B_l$ of $C$ inside $B$ and a~generalized standard family $\mc{C}^{B_l}$ of coray tubes of $\Gamma_{B_l}$ such that $B$ is obtained from $B_l$ (respectively, $\mc{C}^B$ is obtained from $\mc{C}^{B_l}$) by a~sequence of admissible operations of types \adm{1} and \adm{2}, using modules from $\mc{C}^{B_l}$.
Moreover, $\mc{P}^B=\mc{P}^{B^l}$.
Hence $N$ is also a~$B^l$-module and therefore, because $B^l$ is a~quotient algebra of $B$ by an ideal $AhA$ for an idempotent $h$, $M$ is also a~$B^l$-module.
Further, since every component in $\Gamma_{B_l}$ is generalized standard, we infer from \cite{JMS} that $B_l$ is an almost concealed canonical algebra of Euclidean or tubular type.

Assume first that $B_l$ is of Euclidean type.
Then, because every module from $\mc{P}^{B_l}$ is uniquely determined by its composition factors, $N$ belongs to $\mc{C}^{B_l}$ in $\Gamma_{B_l}$.
Thus $N$ is a~module from the family $\mc{C}$ in $\Gamma_A$.

Assume that $B_l$ is of tubular type.
Then $\mc{P}^{B_l}$ consists of all indecomposable modules which precede the family $\mc{C}^{B_l}$ of coray tubes of $\Gamma_{B_l}$.
Hence, applying Proposition \ref{prop_4.1}, we conclude that $N$ belongs to $\mc{C}^{B_l}$.
Thus $N$ is a~module from the family $\mc{C}$ in $\Gamma_A$.

Summing up, the family $\mc{C}^A$ consists of quasi-tubes having common composition factors, is closed under composition factors and, from our assumptions on $\Sigma_A$, consists of modules which do not lie on infinite short cycles. \phantom{.}~\hfill{$\square$}
\end{pf}

Recall the following characterization of self-injective algebras proved in \cite[Theorem 1.1]{KSY}.
\begin{thm}\label{KSY_thm_1.1}
Let $A$ be a~basic, connected, self-injective artin algebra. The
following statements are equivalent.
\begin{enumerate}[(i)]
    \item $\Gamma_A$ admits a~nonempty family $\mc{C}=(\mc{C}_i)_{i\in I}$ of quasitubes having common composition factors, closed on composition factors, and consisting of modules which do not lie on infinite short cycles in $\modd A$.
    \item $A$ is isomorphic to an orbit algebra $\hB/G$, where $B$ is an almost concealed canonical algebra and $G$ is an infinite cyclic group of automorphisms of $\hB$ of one of the forms
    \begin{enumerate}[(a)]
        \item $G=(\varphi\nu_{\hB}^2)$, for a~strictly positive automorphism $\varphi$ of $\hB$,
        \item $G=(\varphi\nu_{\hB}^2)$, for $B$ a~tubular algebra and $\varphi$ a~rigid automorphism of $\hB$,
        \item $G=(\varphi\nu_{\hB}^2)$, for $B$ of Euclidean or wild type and $\varphi$ a~rigid automorphism of $\hB$ acting freely on the nonstable tubes of the unique separating family $\mc{T}^B$ of ray tubes of $\Gamma_B$,
    \end{enumerate}
where $\nu_{\hB}$ is the Nakayama automorphism of $\hB$.
\end{enumerate}
\end{thm}

It follows now from Lemma \ref{lem_4.2} and Theorem \ref{KSY_thm_1.1} that the algebra $A$ is of the form $\hB/(\varphi\nu_{\hB}^2)$, where $B$ is an almost concealed canonical algebra and $\varphi$ is a~positive automorphism of $\hB$.
Moreover, because $\Sigma_A$ contains no short cycles, we infer from \cite{JMS} that $B$ is either a~tilted algebra of Euclidean type or a~tubular algebra.
Thus, in order to prove the Theorem, it remains to show that $\varphi$ is a~strictly positive automorphism of $\hB$ if $B$ is not an exceptional tubular algebra.
This will be a~consequence of Propositions \ref{prop_4.4} and \ref{prop_4.5}.

Note that, by \cite[Corollary 1.4]{KSY}, $k$ is a~field, and therefore $A$ is a~finite-dimensional algebra over a~ field.

We need the following general result which is a~consequence of results proved in \cite{ANS}, \cite{DS}, \cite{Ga}, \cite{NS} and \cite{S1}.
\begin{thm}\label{thm_4.3}
Let $B$ be a~quasi-tilted algebra of canonical type, $G$ an
admissible torsion-free group of automorphisms of $\widehat{B}$,
and $A=\widehat{B}/G$ the associated orbit algebra. Then the
following statements hold.
\begin{enumerate}[(i)]
    \item $G$ is an infinite cyclic group generated by a~strictly positive automorphism $\psi$ of $\widehat{B}$.
    \item The push\dywiz down functor $F_\lambda:\modd\widehat{B}\to\modd A$ associated to the Galois covering $F:\widehat{B}\to\widehat{B}/G=A$ with Galois group $G$ is dense.
    \item The Auslander\dywiz Reiten quiver $\Gamma_A$ of $A$ is isomorphic to the orbit quiver $\Gamma_{\widehat{B}}/G$ of the Auslander\dywiz Reiten quiver $\Gamma_{\widehat{B}}$ of $\widehat{B}$ with respect to the induced action of $G$ on $\Gamma_{\widehat{B}}$.
\end{enumerate}
\end{thm}

\begin{prop}\label{prop_4.4}
Let $B$ be a~non-exceptional tubular algebra, $G$ an infinite cyclic admissible group of automorphisms of $\hB$, and $A=\hB/G$. Then the following statements are equivalent:
\begin{enumerate}[(i)]
    \item The component quiver $\Sigma_A$ of $A$ has no short cycles.
    \item $G=(\varphi\nu^2_{\hB})$, where $\varphi$ is a~strictly positive automorphism of $\hB$.
\end{enumerate}
\end{prop}
\begin{pf}
It follows from the results established in \cite{GeLe}, \cite{HRi2}, \cite{NS} and \cite{S1} (see also \cite{Bi} and \cite{Ku}) that the Auslander\dywiz Reiten quiver $\Gamma_{\hB}$ of $\hB$ has a~decomposition
$$\Gamma_{\hB}=\bigvee_{q\in\bbQ}\mc{C}_q^{\hB}=\bigvee_{q\in\bbQ}\bigvee_{x\in\bbX_q}\mc{C}_{q,x}^{\hB}$$
such that
\begin{enumerate}[(1)]
    \item For each $q\in\bbZ$, $\mc{C}_q^{\hB}$ is an infinite family $\mc{C}_{q,\lambda}^{\hB}$, $\lambda\in\bbX_q$, of quasi-tubes containing at least one projective module.
    \item For each $q\in\bbQ\setminus\bbZ$, $\mc{C}_q^{\hB}$ is an infinite family $\mc{C}_{q,x}^{\hB}$, $x\in\bbX_q$, of stable tubes.
    \item For each $q\in\bbQ$, $\mc{C}_q^{\hB}$ is a~family of pairwise orthogonal generalized standard quasi-tubes with common composition factors, closed under composition factors, and consisting of modules which do not lie on infinite short cycles in $\modd \hB$.
    \item There is a~positive integer $m$ such that $3\leq m\leq\rk K_0(B)$ and $\nu_{\hB}(\mc{C}_q^{\hB})=\mc{C}_{q+m}^{\hB}$ for any $q\in\bbQ$.
    \item $\Hom_{\hB}(\mc{C}_q^{\hB},\mc{C}_r^{\hB})=0$ for all $q>r$ in $\bbQ$.
    \item $\Hom_{\hB}(\mc{C}_q^{\hB},\mc{C}_r^{\hB})=0$ for all $r>q+m$ in $\bbQ$.
    \item For $q\in\bbQ$, we have $\Hom_{\hB}(\mc{C}_q^{\hB},\mc{C}_{q+m}^{\hB})\neq 0$ if and only if $q\in\bbZ$.
    \item For $p<q$ in $\bbQ$ with $\Hom_{\hB}(\mc{C}_p^{\hB},\mc{C}_q^{\hB})\neq 0$, we have $\Hom_{\hB}(\mc{C}_p^{\hB},\mc{C}_r^{\hB})\neq 0$ and $\Hom_{\hB}(\mc{C}_r^{\hB},\mc{C}_q^{\hB})\neq 0$ for any $r\in\bbQ$ with $p\leq r\leq q$.
    \item For all $p\in\bbQ\setminus\bbZ$ and all $q\in\bbQ$ with $\Hom_{\hB}(\mc{C}_p^{\hB},\mc{C}_q^{\hB})\neq 0$, we have $\Hom_{\hB}(\mc{C}_{p,x}^{\hB},\mc{C}_{q,y}^{\hB})\neq 0$ for all $x\in\bbX_p$ and $y\in\bbX_q$.
    \item For all $p\in\bbQ$ and all $q\in\bbQ\setminus\bbZ$ with $\Hom_{\hB}(\mc{C}_p^{\hB},\mc{C}_q^{\hB})\neq 0$, we have $\Hom_{\hB}(\mc{C}_{p,x}^{\hB},\mc{C}_{q,y}^{\hB})\neq 0$ for all $x\in\bbX_p$ and $y\in\bbX_q$.
\end{enumerate}
We know also from Theorem \ref{thm_4.3}~(i) that $G$ is generated by a~strictly positive automorphism $g$ of $\hB$.
Consider the canonical Galois covering $F:\hB\to\hB/G=A$ and the associated push-down functor $F_\lambda:\modd\hB\to\modd A$.
Since $F_\lambda$ is dense, we obtain natural isomorphisms of $k$-modules
$$
\bigoplus_{i\in\bbZ}\Hom_{\hB}(X,{^{g^i}Y})\stackrel{\sim}{\rightarrow}\Hom_A(F_\lambda(X),F_\lambda(Y)),
$$
$$
\bigoplus_{i\in\bbZ}\Hom_{\hB}({^{g^i}X,Y})\stackrel{\sim}{\rightarrow}\Hom_A(F_\lambda(X),F_\lambda(Y)),
$$
for all indecomposable modules $X$ and $Y$ in $\modd\hB$.

We first show that $(ii)$ implies $(i)$.
Assume that $g=\varphi\nu^2_{\hB}$ for some strictly positive automorphism $\varphi$ of $\hB$.
Then it follows from $(4)$ that there is a~positive integer $l> 2m$ such that $g(\mc{C}_q^{\hB})=\mc{C}_{q+l}^{\hB}$ for any $q\in\bbQ$.
Since $g=\varphi\nu_{\hB}^2=(\varphi\nu_{\hB})\nu_{\hB}$ with $\varphi\nu_{\hB}$ a~strictly positive automorphism of ${\hB}$, invoking the knowledge of the supports of indecomposable modules in $\modd {\hB}$ (see \cite[Section 3]{NS}), we conclude that the images $F_\lambda(S)$ and $F_\lambda(T)$ of any non-isomorphic simple ${\hB}$-modules $S$ and $T$ which occur as composition factors of modules in a~fixed family $\mc{C}_q^{\hB}$ are non-isomorphic simple $A$-modules.
Therefore, it follows from Theorem \ref{thm_4.3} and properties $(1)$-$(4)$, that, for each $q\in\bbQ$, $\mc{C}_q^A=F_\lambda(\mc{C}_q^{\hB})=(\mc{C}_{q,x}^A)_{x\in\bbX_q}$, where $\mc{C}_{q,x}^A=F_\lambda(\mc{C}_{q,x}^{\hB})$, $x\in\bbX_q$, is an infinite family of quasi-tubes of $\Gamma_A$ with common composition factors and closed under composition factors.
Take now $p\in\bbQ$.
We claim that $\mc{C}_{p,x}^A$, for any $x\in\bbX_p$, is a~quasi-tube without external short paths in $\modd A$.
Observe first that, for two indecomposable modules $M$ and $N$ in $\mc{C}_p^A$, we have $M=F_\lambda(X)$ and $L=F_\lambda(Y)$, for some indecomposable modules $X$ and $Y$ in $\mc{C}_p^{\hB}$, and $F_\lambda$ induces an isomorphism of $k$-modules $\Hom_A(M,N)\stackrel{\sim}{\rightarrow}\Hom_{\hB}(X,Y)$, by $(5)$, $(6)$ and the inequalities $q+l> q+2m>q+m$.
Suppose now that there is an external short path $M\to L\to N$ in $\modd A$ with $M$ and $N$ in $\mc{C}_{p,x}^A$, for some $x\in\bbX_p$, and $L$ not in $\mc{C}_{p,x}^A$.
Observe that $L$ is not in $\mc{C}_p^A$ because by $(3)$ different quasi-tubes in $\mc{C}_p^A$ are orthogonal.
Therefore, $M=F_\lambda(X)$, $N=F_\lambda(Y)$ for some $X$ and $Y$ in $\mc{C}_{p,x}^{\hB}$ and $L=F_\lambda(Z)$ for some $Z$ in $\mc{C}_r^{\hB}$ with $r>p$.
We have an isomorphism of $k$-modules, induced by $F_\lambda$,
$$
\Hom_A(M,L)\stackrel{\sim}{\rightarrow}\bigoplus_{i\in\bbZ}\Hom_{\hB}(X,{^{g^i}Z}).
$$
Since $\Hom_A(M,L)\neq 0$, we may choose, invoking $(5)$, a~minimal $r>p$ and $Z\in\mc{C}_r^{\hB}$ such that $L=F_\lambda(Z)$ and $\Hom_{\hB}(X,Z)\neq 0$.
Since $X$ lies in $\mc{C}^{\hB}_p$, applying $(6)$ and $(7)$, we infer that $p<r\leq p+m$.
Further, we have also an isomorphism of $k$-modules, induced by $F_\lambda$,
$$
\Hom_A(L,N)\stackrel{\sim}{\rightarrow}\bigoplus_{i\in\bbZ}\Hom_{\hB}(Z,{^{g^i}}Y).
$$
Observe that, for each $i\in\bbZ$, ${^{g^i}Y}$ is an indecomposable module from $\mc{C}_{p+li}^{\hB}$, and clearly with $F_\lambda({^{g^i}Y})=F_\lambda(Y)=N$.
Since $\Hom_A(L,N)\neq 0$, $L=F_\lambda(Z)$ for $Z\in\mc{C}_r^{\hB}$ with $r>p$, and $Y\in\mc{C}_p^{\hB}$, applying $(5)$, we conclude that $\Hom_{\hB}(Z,{^{g^i}Y})\neq 0$, for some $i\geq 1$.
But then $p+li\geq p+l> p+2m\geq r+m$, because $r\leq p+m$, and we obtain a~contradiction with $(6)$.

Summing up, we have proved that all quasi-tubes in $\Gamma_A$ are generalized standard and consists of modules which do not lie on external short paths in $\modd A$.
Thus, the component quiver $\Sigma_A$ of $A$ has no short cycles.
Therefore, $(ii)$ implies $(i)$.

We will show now that that $(i)$ implies $(ii)$.
Assume that the component quiver $\Sigma_A$ has no short cycles.
Then, by Lemma \ref{lem_4.2}, $\Gamma_A$ admits a~family $\mc{C}=(\mc{C}_x)_{x\in\bbX}$ of quasi-tubes with common composition factors, closed under composition factors and consisting of modules which do not lie on infinite short cycles in $\modd A$.
We know from the property $(3)$ that, for each $q\in\bbQ$, $\mc{C}_q^A=F_\lambda(\mc{C}_q^{\hB})$ is a~family $\mc{C}^A_{q,x}=F_\lambda(\mc{C}^{\hB}_{q,x})$, $x\in\bbX_q$, of quasi-tubes with common composition factors.
Moreover, the push-down functor $F_\lambda$ induces an isomorphism of translation quivers $\Gamma_{\hB}/G\stackrel{\sim}{\rightarrow}\Gamma_A$ (see Theorem \ref{thm_4.3}), and hence every component of $\Gamma_A$ is a~quasi-tube of the form $\mc{C}_{q,x}^A=F_\lambda(\mc{C}^{\hB}_{q,x})$ for some $q\in\bbQ$ and $x\in\bbX_q$.
Then, since the family $\mc{C}$ is closed under composition factors, we conclude that there is $r\in\bbQ$ such that $\mc{C}$ contains all quasi-tubes $\mc{C}_{r,x}^A$, $x\in\bbX_r$, of $\mc{C}_r^A$.
This forces, by \cite[Proposition 6.4]{KSY}, $g$ to be of the form $g=\varphi\nu_{\hB}^2$ for some positive automorphism $\varphi$ of ${\hB}$.
Suppose that $\varphi$ is a~rigid automorphism of $\hB$.
Then, from Proposition \ref{prop_3.1}, we know that the restriction of $\varphi$ to $B$ fixes an indecomposable projective module, that is there is an indecomposable projective module $P$ such that $\varphi(P)=P$.
Thus, let $\mc{C}_{p,x}$, for some $p\in\bbZ$ and $x\in\bbX_p$, be the quasi-tube, in $\Gamma_{\hB}$, containing $P$.
Without loss of generality, we may assume that $p=0$.
We have a~short cycle of modules in $\modd\hB$ of the form \begin{tikzcd} P\arrow{r}{f}& \nu_{\hB}(P)\arrow{r}{g}& \nu^2_{\hB}(P)\end{tikzcd}, where $f$ and $g$ are the following compositions of homomorphisms
$$P\to\top(P)\stackrel{\sim}{\rightarrow}\soc(\nu_{\hB}(P))\to\nu_{\hB}(P),$$
and
$$\nu_{\hB}(P)\to\top(\nu_{\hB}(P))\stackrel{\sim}{\rightarrow}\soc(\nu^2_{\hB}(P))\to\nu^2_{\hB}(P).$$
Consequently, we obtain a~short cycle
$$F_\lambda(\mc{C}_{0,x})\to F_\lambda(\mc{C}_{m,y})\to F_\lambda(\mc{C}_{0,x})$$
in $\Sigma_A$, because $\rad_A^\infty(F_\lambda(\mc{C}_{0,x}),F_\lambda(\mc{C}_{m,y}))\neq 0$ and $\rad_A^\infty(F_\lambda(\mc{C}_{m,y}),F_\lambda(\mc{C}_{2m,x}))=\rad_A^\infty(F_\lambda(\mc{C}_{m,y}),F_\lambda(\mc{C}_{0,x}))\neq 0,$ where $\nu_{\hB}(P)\in\mc{C}_{m,\mu}$, for some $y\in\bbX_m$, what contradicts $(i)$. \hfill{$\square$}
\end{pf}

\begin{prop}\label{prop_4.5}
Let $B$ be a~tilted algebra of Euclidean type, $G$ an infinite cyclic admissible group of automorphisms of ${\hB}$, and $A={\hB}/G$. Then the following statements are equivalent:
\begin{enumerate}[(i)]
    \item The component quiver $\Sigma_A$ of $A$ has no short cycle.
    \item $G=(\varphi\nu_{\hB}^2)$, for a~strictly positive automorphism $\varphi$ of ${\hB}$.
\end{enumerate}
\end{prop}
\begin{pf}
It follows from \cite{ANS}, \cite{AS1} and \cite{S1} that the Auslander\dywiz Reiten quiver $\Gamma_{\hB}$ of ${\hB}$ has a~decomposition
$$
\Gamma_{\hB}=\bigvee_{q\in\bbZ}(\mc{C}_q^{\hB}\vee\mc{X}_q^{\hB})
$$
such that
\begin{enumerate}[(1)]
    \item For each $q\in\bbZ$, $\mc{X}_q^{\hB}$ is an acyclic component of Euclidean type.
    \item For each $q\in\bbZ$, $\mc{C}_q^{\hB}$ is a~family $\mc{C}_{q,x}^{\hB}$, $x\in\bbX_q$, of pairwise orthogonal generalized standard quasi-tubes with common composition factors, closed under composition factors, and consisting of modules which do not lie on infinite short cycles in $\modd{\hB}$.
    \item For each $q\in\bbZ$, we have $\nu_{\hB}(\mc{C}_q^{\hB})=\mc{C}_{q+2}^{\hB}$ and $\nu_{\hB}(\mc{X}_q^{\hB})=\mc{X}_{q+2}^{\hB}$.
    \item For each $q\in\bbZ$, we have $\Hom_{\hB}(\mc{X}_q^{\hB},\mc{C}_q^{\hB}\vee\bigvee_{r<q}(\mc{C}_r^{\hB}\vee\mc{X}_r^{\hB}))=0$ and $\Hom_{\hB}(\mc{C}_q^{\hB},\bigvee_{r<q}(\mc{C}_r^{\hB}\vee\mc{X}_r^{\hB}))=0$.
    \item For each $q\in\bbZ$, we have $\Hom_{\hB}(\mc{C}_q^{\hB},\mc{X}_{q+2}^{\hB}\vee\bigvee_{r>q+2}(\mc{C}_r^{\hB}\vee\mc{X}_r^{\hB}))=0$ and $\Hom_{\hB}(\mc{X}_q^{\hB},\bigvee_{r>q+2}(\mc{C}_r^{\hB}\vee\mc{X}_r^{\hB}))=0$.
    \item For $q\in\bbZ$ and $x,y\in\bbX_q$, we have $\Hom_{\hB}(\mc{C}_{q,x}^{\hB},\mc{C}_{q+2,y}^{\hB})\neq 0$ if and only if the quasi-tube $\mc{C}_{q,x}^{\hB}$ is non-stable and $\nu_{\hB}(\mc{C}_{q,y}^{\hB})=\mc{C}_{q+2,x}^{\hB}$.
    \item For all $q\in\bbZ$ and $x,y\in\bbX_q$, we have $\Hom_{\hB}(\mc{C}_{q,x}^{\hB},\mc{C}_{q+1,y}^{\hB})\neq 0$.
    \item Each each $r\in\bbZ$, $\mc{X}_r$ contains at least one projective module.
\end{enumerate}

We know also from Theorem \ref{thm_4.3} that $G$ is generated by a~strictly positive automorphism $g$ of ${\hB}$.
Hence there exists a~positive integer $l$ such that $g(\mc{C}_q^{\hB})=\mc{C}_{q+l}^{\hB}$ and $g(\mc{X}_q^{\hB})=\mc{X}_{q+l}^{\hB}$ for any $q\in\bbZ$.
Consider the canonical Galois covering $F:{\hB}\to{\hB}/G=A$ and the associated push-down functor $F_\lambda:\modd{\hB}\to\modd A$.
Since $F_\lambda$ is dense, we obtain natural isomorphisms of
$k$-modules
$$
\bigoplus_{i\in\bbZ}\Hom_{\hB}(X,{^{g^i}Y})\stackrel{\sim}{\rightarrow}\Hom_A(F_\lambda(X),F_\lambda(Y)),
$$
$$
\bigoplus_{i\in\bbZ}\Hom_{\hB}({^{g^i}X},Y)\stackrel{\sim}{\rightarrow}\Hom_A(F_\lambda(X),F_\lambda(Y)),
$$
for all indecomposable modules $X$ and $Y$ in $\modd{\hB}$.

We show first that $(i)\Rightarrow(ii)$.
Assume that the component quiver $\Sigma_A$ has no short cycles.
Then, by Lemma \ref{lem_4.2}, $\Gamma_A$ admits a~family $\mc{C}=(\mc{C}_\lambda)_{\lambda\in\bbX}$ of quasi-tubes with common composition factors, closed under composition factors and consisting of modules which do not lie on infinite short cycles in $\modd A$.
Then it follows from \cite[Proposition 6.5]{KSY} that $g=\varphi\nu^2_{\hB}$ for some positive automorphism $\varphi$ of $\hB$.
We claim that $\varphi$ is a~strictly positive automorphism of $\hB$.

Assume that $\varphi$ is a~rigid automorphism of ${\hB}$.
Take $q=0$ and, invoking $(8)$, some projective-injective module $P$ in $\mc{X}_0$.
Let $f$ and $g$ be the following compositions of homomorphisms
$$P\to\top(P)\stackrel{\sim}{\rightarrow}\soc(\nu_{\hB}(P))\to\nu_{\hB}(P)$$
and
$$\nu_{\hB}(P)\to\top(\nu_{\hB}(P))\stackrel{\sim}{\rightarrow}\soc(\nu^2_{\hB}(P))\to\nu^2_{\hB}(P),$$
respectively.
Then we have a~short path of indecomposable modules
\begin{center}
\begin{tikzcd} P\arrow{r}{f} & \nu_{\hB}(P)\arrow{r}{g} &\nu^2_{\hB}(P)\end{tikzcd}
\end{center}
in $\modd\hB$, where, by $(3)$, $P\in\mc{X}_0$, $\nu_{\hB}(P)\in\mc{X}_2$ and $\nu^2_{\hB}(P)\in\mc{X}_4$.
Thus, it follows, from Theorem \ref{thm_4.3}, that we have a~short path of indecomposable modules $F_\lambda(P)\to F_\lambda(\nu_{\hB}(P))\to F_\lambda(\nu^2_{\hB}(P))$ in $\modd A$.
Because $\varphi$ is a~rigid automorphism of $\hB$ we conclude that $F_\lambda(P)$ and $F_\lambda(\nu^2_{\hB}(P))$ belong to the same component $F_\lambda(\mc{X}_0)$.
Obviously, $\rad^\infty_A(F_\lambda(P),F_\lambda(\nu_{\hB}(P))\neq 0$ and $\rad^\infty_A(F_\lambda(\nu_{\hB}(P),F_\lambda(\nu^2_{\hB}(P))\neq 0$.
Therefore, the component quiver $\Sigma_A$ contains a~short cycle $F_\lambda(\mc{X}_0)\to F_\lambda(\mc{X}_2)\to F_\lambda(\mc{X}_0)$, and we get a~contradiction.

This finishes the proof that $(i)$ implies $(ii)$.

Assume now that $(ii)$ holds.
In particular, we have $g=\varphi\nu_{\hB}^2$ for a~strictly positive automorphism $\varphi$ of ${\hB}$.
Then it follows from $(3)$ that there is a~positive integer $l>4$ such that $g(\mc{C}_q^{\hB})=\mc{C}_{q+l}^{\hB}$ and $g(\mc{X}_q^{\hB})=\mc{X}_{q+l}^{\hB}$ for any $q\in\bbZ$.
By $(2)$ and Theorem \ref{thm_4.3}, in order to show that $\Sigma_A$ has no short cycles, we must show, that every component in $\Gamma_A$ is generalized standard and has no external short paths.
From property $(2)$ and \cite[Theorem 3]{SY3} every component in $\Gamma_A$ is generalized standard.
Assume that there is a~component $\mc{C}$ in $\Gamma_A$ and an external short path $M\to N\to L$ with $M$ and $L$ in $\mc{C}$ but $N$ not in $\mc{C}$.
By Theorem \ref{thm_4.3}, there are indecomposable modules $X,Y$ and $Z$ in $\modd\hB$ such that $M=F_\lambda(X)$, $N=F_\lambda(Y)$ and $L=F_\lambda(Z)$.
Moreover, $X$ belongs to $\mc{C}_{p,x}$, for some $p\in\bbZ$ and $x\in\bbX_p$, or $X$ belongs to $\mc{X}_p$, for some $p\in\bbZ$.
Then $\mc{C}=F_\lambda(\mc{C}_{p,x})$ or $\mc{C}=F_\lambda(\mc{X}_p)$.
Without loss of generality, we may assume that $p=0$.
Therefore, we have two cases to consider.

Assume that $X\in\mc{C}_{0,x}$.
We have an isomorphism of $k$-modules, induced by $F_\lambda$, $$\Hom_A(M,N)\stackrel{\sim}{\rightarrow}\bigoplus_{i\in\bbZ}\Hom_{\hB}(X,{^{g^i}Y}).$$
Since $\Hom_A(M,N)\neq 0$, invoking $(5)$, we conclude that $Y$ belongs to $$\mc{X}_0\vee\mc{C}_{1}\vee\mc{X}_{1}\vee\mc{C}_{2}.$$
We have also an isomorphism of $k$-modules $$\Hom_A(N,L)\stackrel{\sim}{\rightarrow}\bigoplus_{i\in\bbZ}\Hom_{\hB}(Y,{^{g^i}Z}).$$
Again, since $\Hom_A(N,L)\neq 0$, we conclude from $(5)$ that $Z$ belongs to $$\mc{X}_2\vee\mc{C}_{2}\vee\mc{X}_3\vee\mc{C}_{3}\vee\mc{X}_4\vee\mc{C}_{4}.$$
On the other hand, by the property $(4)$ and our assumption on $\varphi$, we have that $Z\in\mc{C}_{l,\lambda}$, for some $l>4$, a~contradiction.

Assume that $X\in\mc{X}_0$.
We have an isomorphism of $k$-modules, induced by $F_\lambda$, $$\Hom_A(M,N)\stackrel{\sim}{\rightarrow}\bigoplus_{i\in\bbZ}\Hom_{\hB}(X,{^{g^i}Y}).$$
Since $\Hom_A(M,N)\neq 0$, invoking $(5)$, we infer that $Y$ belongs to $$\mc{C}_1\vee\mc{X}_{1}\vee\mc{C}_{2}\vee\mc{X}_{2}.$$
We have also an isomorphism of $k$-modules, $$\Hom_A(N,L)\stackrel{\sim}{\rightarrow}\bigoplus_{i\in\bbZ}\Hom_{\hB}(Y,{^{g^i}Z}).$$
Again, since $\Hom_A(N,L)\neq 0$, we conclude from $(5)$ that $Z$ belongs to $$\mc{X}_{1}\vee\mc{X}_{2}\vee\mc{X}_{3}\vee\mc{X}_{4}.$$
On the other hand, by the property $(4)$ and our assumption on $\varphi$, we obtain that $Z\in\mc{X}_l$, for some $l>4$, a~contradiction. $\phantom{.}$~{\hfill$\square$}
\end{pf}

\section*{Acknowledgements}

The project was supported by the Polish National Science Center grant awarded on the basis of the decision number DEC-2011/01/N/ST1/02064.

\bibliographystyle{mk}
\bibliography{pub_3}
\end{document}